\documentclass[a4paper,11pt]{amsart}
\usepackage[arrow,matrix]{xy}
\usepackage{amsmath,amssymb,bbm, amscd, amsthm,mathrsfs}
\usepackage[colorlinks=true,citecolor=blue,linkcolor=blue]{hyperref}
\theoremstyle{plain}
\newtheorem{thm}{Theorem}[section]
\newtheorem{cor}[thm]{Corollary}
\newtheorem{lem}[thm]{Lemma}
\newtheorem{prop}[thm]{Proposition}

\newtheorem{exa}[thm]{Example}

\newtheorem{rem}[thm]{Remark}

\def\im{\operatorname {im}}    \def\ot{\otimes}
\def\Hom{\operatorname {Hom}}
\def\RHom{\operatorname {RHom}}
\def\Ext{\operatorname {Ext}}
\def\Tor{\operatorname {Tor}}\def\g{\mathfrak{g}}
\def\k{\mathrm{I\hspace{-2pt}k}}

\begin{document}
\title{\bf Cocommutative Calabi-Yau Hopf algebras and deformations}

\author{Ji-Wei He, Fred Van Oystaeyen and Yinhuo Zhang}
\address{J.-W. He\newline \indent Department of Mathematics, Shaoxing College of Arts and Sciences, Shaoxing Zhejiang 312000,
China\newline \indent Department of Mathematics and Computer
Science, University of Antwerp, Middelheimlaan 1, B-2020 Antwerp,
Belgium} \email{jwhe@zscas.edu.cn}
\address{F. Van Oystaeyen\newline\indent Department of Mathematics and Computer
Science, University of Antwerp, Middelheimlaan 1, B-2020 Antwerp,
Belgium} \email{fred.vanoystaeyen@ua.ac.be}
\address{Y. Zhang\newline
\indent Department WNI, University of Hasselt, Universitaire Campus,
3590 Diepenbeek, Belgium} \email{yinhuo.zhang@uhasselt.be}

\date{}

\begin{abstract} The Calabi-Yau property of cocommutative Hopf
algebras is discussed by using the homological integral, a recently
introduced tool for studying infinite dimensional AS-Gorenstein Hopf
algebras. It is shown that the skew-group algebra of a universal
enveloping algebra of a finite dimensional Lie algebra $\g$ with a
finite subgroup $G$ of automorphisms of $\g$ is Calabi-Yau if and
only if the universal enveloping algebra itself is Calabi-Yau and
$G$ is a subgroup of the special linear group $SL(\g)$. The
Noetherian cocommutative Calabi-Yau Hopf algebras of dimension not
larger than 3 are described. The Calabi-Yau property of Sridharan
enveloping algebras of finite dimensional Lie algebras is also
discussed. We obtain some equivalent conditions for a Sridharan
enveloping algebra to be Calabi-Yau, and then partly answer a
question proposed by Berger. We list all the nonisomorphic
3-dimensional Calabi-Yau Sridharan enveloping algebras.

\end{abstract}


\maketitle

\section*{Introduction}

We work over a fixed field $\k$ which is assumed to be
algebraically closed and of characteristic zero and is assumed to be
the field of complex numbers $\mathbb{C}$ if necessary.

Calabi-Yau algebras are studied in recent years because of their
applications in algebraic geometry and mathematical physics. The aim
of this paper is to try to understand cocommutative Calabi-Yau Hopf
algebras of lower dimensions. We take the Calabi-Yau (CY) property
from Ginzburg \cite{G}, the definition will be recalled in Section
\ref{sec2}. The main tool used in this paper is the homological
integral of an AS-Gorenstein Hopf algebra recently introduced by
Lu-Wu-Zhang in \cite{LWZ} and extended to general AS-Gorenstein
algebras (the definition is recalled in Section \ref{sec1}) by
Brown-Zhang in \cite{BZ}. As a consequence of Kostant-Larson's work
(cf. \cite{La,Sw}), we know that any cocommutative Hopf algebra
(over an algebraically closed field) is isomorphic to a skew-group
(or smash product) algebra of the universal enveloping subalgebra of
the primitive elements and the group subalgebra of the group-like
elements. Hence it is necessary to discuss how the homological
integral works on the skew-group algebras. In Section \ref{sec1}, we
discuss finite group actions on AS-Gorenstein algebras. Let $A$ be
an AS-Gorenstein algebra and $G$ be a finite group. If there is a
$G$-action on $A$ such that the $G$-action is compatible with the
augmentation map of $A$, then $A$ is called an augmented $G$-module
algebra. If $A$ is an augmented $G$-module algebra which is
AS-Gorenstein, then we show in Section \ref{sec1} that the the
skew-group algebra $A\#\k G$ is also AS-Gorenstein and the (left)
homological integral $\int^l_A$ is a left $G$-module and the (left)
homological integral of $A\#\k G$ is equal to the $G$-invariants of
the left $G$-module $\int^l_A\#\k G$ (Proposition \ref{prop1}). A
necessary and sufficient condition for a Noetherian Hopf algebra to
be CY is given in Section \ref{sec2}. It turns out that a Noetherian
CY Hopf algebra must be AS-regular and has bijective antipode. When
the AS-Gorenstein Hopf algebra is the universal enveloping algebra
of a finite dimensional Lie algebra, we have the following result
(Theorem \ref{thm1}).

\noindent{\bf Theorem}. {\it Let $\mathfrak{g}$ be a finite
dimensional Lie algebra, and $G\subseteq Aut_{Lie}(\mathfrak{g})$ be
a finite group. Then the skew-group algebra
$U(\mathfrak{g})\#\k G$ is CY of dimension $d$ if and only
if $U(\mathfrak{g})$ is CY of dimension $d$ and $G\subseteq
SL(\mathfrak{g})$.}

Let $H$ be a Hopf algebra, $G(H)$ the group of the group-like
elements of $H$ and $P(H)$ the space of the primitive elements of
$H$. Applying the above theorem, we can list all cocommutative CY
Hopf algebras $H$ of global dimension not larger than 3 such that
$G(H)$ is finite and $P(H)$ is finite dimensional. It is easy to
determine the 1-dimensional Noetherian CY cocommutative Hopf
algebras with finite group $G(H)$. They are the tensor products of
the polynomial algebra $\k [x]$ with the group algebras of some
finite groups. In the 2-dimensional case, such a Hopf algebras $H$
must be isomorphic to a skew-group algebra of the form $H\cong\k
[x,y]\#\k G$, where $G$ is a finite group and the $G$-action on
$\k[x,y]$ is induced by a group map $\nu:G\to SL(2,\k )$ (Theorem
\ref{thm2}).

For the 3-dimensional case, we show that there are only 4 cases of
nonisomorphic 3-dimensional Lie algebras whose universal enveloping
algebras are CY (Proposition \ref{prop5}).  A 3-dimensional
Noetherian CY cocommutative Hopf algebra $H$ with finite group
$G(H)$ and finite dimensional $P(H)$ is isomorphic to a skew-group
algebra of form $U(\g)\# \k G$, where the Lie algebra $\g$ is one of
the Lie algebras listed in Proposition \ref{prop5} and $G$ is a
finite group with a group morphism $\nu:G\to Aut_{Lie}(\g)$ such
that $\im(\nu)$ is  a subgroup of $SL(\g)$ (Theorem \ref{thm5}).

In the last section, we deal with the Sridharan enveloping algebras
of finite dimensional Lie algebras. In general, a Sridharan
enveloping algebra is no longer a Hopf algebra. However, a Sridharan
enveloping algebra is a cocycle deformation of a cocommutative Hopf
algebras, and the CY property of the Sridharan enveloping algebras
is closely related to that of the universal enveloping algebras.
Hence it is proper to include the discussion of Sridharan enveloping
algebras in this paper. Sridharan enveloping algebras were
introduced in \cite{Sr} in order to discuss certain representations
of Lie algebras. Let $\g$ be a finite dimensional Lie algebra. A
Sridharan enveloping algebra is related to a 2-cocycle $f\in
Z^2(\g,\k)$ of $\g$, and is usually denoted by $U_f(\g)$ (the
definition is recalled in the final section). The class of Sridharan
enveloping algebras includes many interesting algebras, such as Weyl
algebras. Homological properties, especially the Hochschild
(co)homology and cyclic homology, are studied by several authors
\cite{Sr,Ka,N0,N}. Berger proposed at the end of his recent paper
\cite{Be} a question: to find some necessary and sufficient
conditions for a Sridharan enveloping algebra to be CY. We get the
following result (Theorem \ref{thm3}) which in part answers Berger's
question.

\noindent{\bf Theorem}. {\it Let $\g$ be a finite dimensional Lie
algebra, and $f\in Z^2(\g,\k)$ be an arbitrary 2-cocycle of $\g$.
The following statements are equivalent.
\begin{itemize}
\item [(i)] The Sridharan
enveloping algebra $U_f(\g)$ is CY of dimension $d$.
\item [(ii)] The universal enveloping algebra $U(\g)$ is CY of dimension
  $d$.
\item [(iii)] $\dim\g=d$, and for any $x\in\g$, $\text{\rm tr}(\text{\rm ad}_\g(x))=0$.
\end{itemize}}

We also list all the 3-dimensional CY Sridharan enveloping algebras
at the end of the paper (Theorem \ref{thm4}). There are exactly 7
classes of nonisomorphic 3-dimensional CY Sridharan enveloping
algebras.

\section{Homological integrals of skew group algebras}\label{sec1}

Let $A$ be a left Noetherian augmented algebra with a fixed
augmentation map $\varepsilon:A\longrightarrow \k $. Recall
that $A$ is said to be {\it left Artin-Schelter Gorenstein}
(AS-Gorenstein for short, cf. \cite{BZ}), if
\begin{itemize}
   \item [(i)] injdim${}_AA=d<\infty$,
   \item [(ii)] $\dim\Ext^d_A({}_A\k ,{}_AA)=1$ and
   $\Ext^i_A({}_A\k ,{}_AA)=0$ for all $i\neq d$.
\end{itemize}

If $A$ is a right Noetherian augmented algebra with a fixed
augmentation map, and the right versions of (i) and (ii) above hold,
then $A$ is said to be {\it right AS-Gorenstein}. If $A$ is both
left and right AS-Gorenstein (relative to the same augmentation
map $\varepsilon$), then we say that $A$ is {\it AS-Gorenstein}. Furthermore,
if gldim$A<\infty$, then $A$ is called an {\it AS-regular} algebra.

The concept of a homological integral was first introduced in
\cite{LWZ} for an AS-Gorenstein Hopf algebra as a generalization of
the classical concept of an integral for a finite dimensional Hopf
algebra. The concept was further extended to a general AS-Gorenstein
algebra in \cite{BZ}. It seems that the homological integral is a
useful tool to study infinite dimensional noncocommutative algebras.
Let $A$ be a left AS-Gorenstein algebra. Then $\Ext^d_A({}_A\k
,{}_AA)$ is a one-dimensional right $A$-module. Any nonzero element
in $\Ext^d_A({}_A\k ,{}_AA)$ is called a {\it left homological
integral} of $A$. Write $\int^l_A$ for $\Ext^d_A({}_A\k ,{}_AA)$,
and call it, by a slight abuse of terminology, {\it the homological
integral of $A$}. Similarly, if $A$ is a right AS-Gorenstein
algebra, any nonzero element of the one-dimensional left module
$\Ext^d_A(\k _A,A_A)$ is called a {\it right homological integral}
of $A$. We donote it by $\int_A^r$.

Let $G$ be a finite group. A left $G$-module algebra is called {\it
an augmented $G$-module algebra} if $A$ has an augmentation map
$\varepsilon$ and $\varepsilon$ is a $G$-map. For an augmented left
$G$-module algebra $A$, the skew group algebra $A\#\k G$ is also an
augmented algebra with the augmentation map $\epsilon:A\#\k
G\longrightarrow \k $ defined by $a\#g\mapsto \varepsilon(a)$ for
all $a\in A$ and $g\in G$. With this augmentation map, $\k $ is
naturally a $A\#\k G$-$A\#\k G$-bimodule. Since $G$ is a finite
group, $A\#\k G$ is a left Noetherian algebra if $A$ is left
Noetherian.

Let $M$ and $N$ be left $A\#\k G$-modules. Then $\Hom_A(M,N)$ is a left $G$-module with the adjoint action:
$$g\rightharpoonup f(m)=g\cdot f(g^{-1}\cdot m),$$ for $g\in G$,
$f\in\Hom_A(M,N)$ and $m\in M$. For a left $G$-module $X$, let
$X^G=\{x\in X|g\cdot x=x,\ \text{for all $g\in G$}\}$ be the set
of $G$-invariant elements. It is easy to see that
\begin{equation}\label{tag1}
    \Hom_{A\#\k G}(M,N)=\Hom_A(M,N)^G.
\end{equation}
Since $G$ is finite, the functor $(-)^G$ is exact. It follows that a
left $A\#\k G$-module $Q$ is projective (resp. injective) if and
only if it is projective (resp. injective) as an $A$-module. Also
the $G$-module structure on $\Hom_A(M,N)$ can be extended to the
extension groups $\Ext^i_A(M,N)$, and the isomorphism (\ref{tag1})
can be extended to the following isomorphisms (cf. \cite{MMM})
\begin{equation}\label{tag0}
    \Ext^i_{A\#\k G}(M,N)\cong\Ext^i_A(M,N)^G,\qquad \text{for all
$i\ge0$}.
\end{equation}

Now let $A$ be a left AS-Gorenstein algebra of injdim${}_AA=d$
 and consider the one dimensional module
${}_{A\#\k G}\k $ with the module structure defined
by the augmentation map $\epsilon$. Let
$$\cdots\longrightarrow P^{-n}\overset{\partial^{-n}}\longrightarrow\cdots\overset{\partial^{-2}}\longrightarrow
P^{-1} \overset{\partial^{-1}}\longrightarrow P^0\longrightarrow
{}_{A\#\k G}\k \longrightarrow0$$ be a finitely
generated projective resolution of the $A\#\k G$-module
${}_{A\#\k G}\k $. Then the resolution can be regarded as a
projective resolution of the $A$-module ${}_A\k $. Applying
the functor $\Hom_A(-,A\#\k G)$ to the projective resolution
above, we obtain a complex
\begin{equation}\label{tag2}
\begin{array}{cl}
    \cdots\longleftarrow
\Hom_A(P^{-n},A\#\k G)\longleftarrow \cdots & \longleftarrow
\Hom_A(P^{-1},A\#\k G) \\
& \longleftarrow  \Hom_A(P^0,A\#\k G)\longleftarrow 0
\end{array}
\end{equation}
Since $G$ is a finite group, there are natural isomorphisms of vector
spaces, for $n\ge0$,
\begin{equation}\label{tag3}
    \varphi^n:\Hom_A(P^{-n},A)\ot \k G\longrightarrow\Hom_A(P^{-n},A\#\k G)
\end{equation}
defined by $\varphi^n(f\ot g)(p)=f(p)\#g$ for
$f\in\Hom_A(P^{-n},A)$, $g\in G$ and $p\in P^{-n}$.

Let $Y_A$ be an $A$-module. The tensor space $Y\ot \k G$ is a right
$A\#\k G$-module defined by $$(y\ot g)\cdot(a\#h)=y\cdot(ga)\ot gh,\
\text{ for $y\in Y$, $g,h\in G$ and $a\in A$.}$$ We write this right
$A\#\k G$-module as $Y\#\k G$. Since $\Hom_A(P^{-n},A)$ is a right
$A$-module, $\Hom_A(P^{-n},A)\ot \k G$ is a right $A\#\k G$-module.
$\Hom_A(P^{-n},A\#\k G)$ is also a right $A\#\k G$-module. Thus the
natural isomorphisms in (\ref{tag3}) are in fact right $A\#\k
G$-module isomorphisms
\begin{equation}\label{tag4}
\varphi^n:\Hom_A(P^{-n},A)\#
\k G\longrightarrow\Hom_A(P^{-n},A\#\k G).
\end{equation}
Observe that $A$ is a left $A\#\k G$-module and $\Hom_A(P^{-n},A)$
is a left $G$-module. With the diagonal $G$-action,
$\Hom_A(P^{-n},A)\#\k G$ becomes a left $G$-module. On the other
side, $\Hom_A(P^{-n},A\#\k G)$ is also a left $G$-module with the
adjoint $G$-action. Thus it is not hard to see that both
$\Hom_A(P^{-n},A)\#\k G$ and $\Hom_A(P^{-n},A\#\k G)$ are left $G$-
and right $A\#\k G$-bimodules, and the isomorphisms $\varphi^n$ in
(\ref{tag4}) are isomorphisms of $G$-$A\#\k G$-bimodules.

Now one may check that the complex (\ref{tag2}) is a complex of
$G$-$A\#\k G$-bimodules, and it is isomorphic to the following
complex of $G$-$A\#\k G$-bimodules {\small$$\cdots\longleftarrow
\Hom_A(P^{-n},A)\#\k G\overset{(\partial^{-n})^*\ot
id}\longleftarrow\cdots\overset{(\partial^{-2})^*\ot
id}\longleftarrow \Hom_A(P^{-1},A)\#\k G$$
$$\hspace{70mm}\overset{(\partial^{-1})^*\ot id}\longleftarrow
\Hom_A(P^0,A)\#\k G\longleftarrow0.$$} By taking the
cohomologies of the complex (\ref{tag2}) and those of the complex
above we arrive at isomorphisms of $G$-$A\#\k G$-bimodules:
\begin{equation}\label{tag5}
\Ext^i_A({}_A\k ,A\#\k G)\cong\Ext^i_A({}_A\k ,{}_AA)\#\k G
\end{equation}
for all $i\ge0$. Hence by (\ref{tag0}), we have right
$A\#\k G$-module isomorphisms
\begin{equation}\label{tag6}
\begin{array}{ccl}
\Ext_{A\#\k G}^i({}_{A\#\k G}\k ,A\#\k G)& \cong& \Ext^i_A({}_A\k ,A\#\k G)^G\\
 & \cong & (\Ext^i_A({}_A\k ,{}_AA)\#\k G)^G
 \end{array}
\end{equation}
for all $i\ge0$.

Summarizing the above we arrive at the following results.

\begin{prop}\label{prop1} Let $G$ be a finite group, $A$ an augmented left $G$-module
algebra. Assume $A$ is a left AS-Gorenstein algebra with {\rm
injdim}${}_AA=d$. Then the following statements hold.
\begin{enumerate}
\item[(i)] The left homological integral $\int_A^l$ is a
one-dimensional left $G$-module, and the $G$-action is compatible
with the right $A$-module structure of $\int^l_A$, that is; for
$g\in G$, $a\in A$ and $t\in \int^l_A$, $g(ta)=(gt)(ga)$.
\item[(ii)] $A\#kG$ is left AS-Gorenstein of \text{\rm
injdim}${}_{A\#\k G}(A\#\k G)=d$, and
 as right $A\#\k G$-modules\\
\centerline{$\int^l_{A\#\k G}\cong(\int^l_A\#\k G)^G,$}\\
where the $G$ acts on $\int^l_A\#\k G$ diagonally.
\end{enumerate}
\end{prop}
\proof The statement (i) is evident, and the isomorphism in (ii) is
a direct consequence of the isomorphisms in (\ref{tag6}) if $A\#\k
G$ is left AS-Gorenstein.

What remains to be shown is that the left injective dimension of
$A\#\k G$ is $d$, \linebreak and
$\dim\Ext_{A\#\k G}^d({}_{A\#\k G}\k ,A\#\k G)=1$.
Let
\begin{equation}\label{tag7}
    0\longrightarrow A\longrightarrow
Q^0\overset{\delta^0}\longrightarrow
Q^1\overset{\delta^1}\longrightarrow\cdots
\overset{\delta^{d-1}}\longrightarrow
Q^d\overset{\delta^d}\longrightarrow\cdots
\end{equation}
be an injective resolution of  ${}_{A\#\k G}A$. Since $G$ is
a finite group, all the $Q^i$'s are injective as left $A$-modules.
Hence $coker\delta^{d-1}$ is injective as an $A$-module by the
assumption that ${}_AA$ has injective dimension $d$, which in turn
implies that $coker\delta^{d-1}$ is injective as an
$A\#\k G$-module. Thus we may assume that the resolution
(\ref{tag7}) ends at the $d$-th position. Now tensoring (\ref{tag7}) with
$\k G$, we obtain an exact sequence
$$0\longrightarrow A\ot
\k G\longrightarrow Q^0\ot \k G\longrightarrow
Q^1\ot \k G\longrightarrow\cdots \longrightarrow Q^d\ot
\k G\longrightarrow 0. $$
For a left $A\#\k G$-module $M$, the space $M\ot\k G$ is a left $A\#\k G$-module defined by\\
\centerline{ $(a\#g)\cdot (m\ot h)=a(gm)\ot gh$.} \\
Since $Q^i$ is injective as an $A$-module, $Q^i\ot\k G$ is injective
as an $A\#\k G$-module for all $i\ge0$. Therefore we obtain that the
injective dimension of ${}_{A\#\k G}A\#\k G$ is not larger than $d$.
We claim that $\Ext_{A\#\k G}^d({}_{A\#\k G}\k ,A\#\k
G)\cong(\int^l_A\#\k G)^G\neq0$. Assume that $\alpha$ is a nonzero
element in $\int^l_A$. Since $\dim\int^l_A=1$, there is an algebra
map $\pi: \k G\longrightarrow \k $ such that
$g\cdot\alpha=\pi(g)\alpha$ for all $g\in G$. Let $\pi^{-1}$ be the
convolution inverse of $\pi$ in the dual Hopf algebra $\k G^*$. Then
$\pi^{-1}$ is an algebra map from $\k G$ to $\k $. Hence $\pi^{-1}$
defines an one-dimensional $G$-module. Since $G$ is a finite group,
we have that there is an element $0\neq t\in \k G$ such that
$gt=\pi^{-1}(g)t$ for all $g\in G$. Now for $g\in G$,
$g\cdot(\alpha\#t)=(g\cdot\alpha)\#gt=\pi(g)\pi^{-1}(g)\alpha\#t=\alpha\#t$.
The claim follows. Therefore injdim${}_{A\#\k G}A\#\k G=d$. If
$0\neq t'\in \k G$ is another element such that $g\cdot
(\alpha\#t')=\alpha\#t'$ for all $g\in G$, then we get
$gt'=\pi(g)^{-1}t'=\pi^{-1}(g)t'$ for all $g\in G$. Then we must
have $t'=kt$ for some $k\in\k $. Otherwise, there would be two
one-dimensional $G$-modules that are isomorphic to each other in the
decomposition of the regular $G$-module, which certainly contradicts
the well-known result that in the decomposition of the regular
representation of a finite group the multiplicity of an irreducible
representation equals its dimension (cf. \cite[Sec. 2.4]{S}). Hence
$\dim\Ext_{A\#\k G}^d({}_{A\#\k G}\k ,A\#\k G)=1$. \qed

We certainly should not expect that $\int^l_{A\#\k G}\cong\k\alpha\#
t$ such that $\alpha\in\int^l_A$ and $t$ is an integral of $\k G$,
see the following example.

\begin{exa} {\rm Let $G$ be a cyclic group of order $p$. Let
$\lambda$ be a generator of $G$. Assume that $A$ is an augmented
$G$-module algebra and $A$ is left AS-Gorenstein. Then $A\#\k G$ is
left AS-Gorenstein. In fact, if $G$  acts on $\int_A^l$ trivially,
that is, $\lambda\cdot \alpha=\alpha$ where
$0\neq\alpha\in\int_A^l$, then $\Ext^d_{A\#\k G}({}_{A\#\k G}\k
,A\#\k G)=\k \alpha\# t$ where
$t=\frac{1}{p}\sum_{i=0}^{p-1}\lambda^i$. Now suppose $G$  acts on
$\int^l_A$ nontrivially. Then $\lambda\cdot\alpha=\omega\alpha$,
where $\omega\in \k$ is a $p$th root of the unit. Assume that
$t'=x_0\lambda^0+x_1\lambda+\cdot+x_{p-1}\lambda^{p-1}$ is such that
$\lambda\cdot(\alpha\#t')=\alpha\#t'$. Then we obtain
$$\left\{\begin{array}{lcl}
      x_0 &=& \omega x_{p-1} \\
     x_1 &=& \omega x_0 \\
         &\vdots & \\
     x_{p-1} &=& \omega x_{p-2}.
     \end{array}\right.$$ Obviously, the linear equations above have
      an 1-dimensional solution space with a basis given by $(x_0,x_1,\dots,x_{p-2},x_{p-1})=(\omega,\omega^2,\dots,\omega^{p-1},1)$.
Let
$t'=\omega\lambda^0+\omega^2\lambda+\cdots+\omega^{p-1}\lambda^{p-2}+\lambda^{p-1}$.
Then one can check that\\
\centerline{$\Ext^d_{A\#\k G}({}_{A\#\k G}\k ,A\#\k G)\cong(\int_A^l\#\k G)^G=\k \alpha\#t'$.}
}
\end{exa}

Let $A$ be an augmented left $G$-module algebra. If $A$ is right
AS-Gorenstein, we want to know what $\int^r_{A\#\k G}$ looks like.
The right version of Proposition \ref{prop1} also holds, but it is
more complicated. The algebra $A$ can be viewed as an augmented
right $G$-module algebra through the right $G$-action: $a\cdot
g=g^{-1}a$ for $a\in A$ and $g\in G$. We have the skew group algebra
$\k G\#A$ defined in the usual way. There is an algebra isomorphism
$\theta:A\#\k G\longrightarrow \k G\#A$ by $a\#g\mapsto g\#g^{-1}a$.
Moreover, $\theta$ is compatible with the augmentation maps of
$A\#\k G$ and $\k G\#A$ respectively. Now we can deal with right
$A\#\k G$-modules as right $\k G\#A$-modules.  Let $M$ and $N$ be
right $\k G\# A$-modules. $\Hom_{A}(M,N)$ is a right $G$-module
through the $G$-action: $(f\leftharpoonup g)(m)=f(mg^{-1})g$ for
$f\in\Hom_{A}(M,N)$, $g\in G$ and $m\in M$. Also we have $\Hom_{\k
G\#A}(M,N)=\Hom_{A}(M,N)^G$. Similar to Proposition \ref{prop1}, we
have:

\begin{prop}\label{prop2} Let $G$ be a finite group, $A$ be an augmented left $G$-module
algebra. Assume $A$ is a right AS-Gorenstein algebra with {\rm
injdim}$A_A=d$. Then the following statements hold.
\begin{enumerate}
\item[(i)] The right homological integral $\int_A^r$ is a
1-dimensional right $G$-module, and the $G$-action is compatible
with the left $A$-module structure of $\int^r_A$, that is; for $g\in
G$, $a\in A$ and $t\in \int^r_A$, $(at)\cdot g=(a\cdot g)(t\cdot
g)=(g^{-1}a)(t\cdot g)$.
\item[(ii)] $A\#kG$ is right AS-Gorenstein and \text{\rm
injdim}$(A\#\k G)_{A\#\k G}=d$, also as left
$A\#\k G$-modules:\\
\centerline{$\int^r_{A\#\k G}\cong(\k G\ot\int^r_A)^G,$}\\
where the left $A\#\k G$-action on $\k G\ot\int^r_A$
is given by $(a\#g)\cdot(h\ot \alpha)=gh\ot (h^{-1}g^{-1}a)\alpha$
for $g,h\in G$, $a\in A, \alpha\in\int^r_A$, and the right
$G$-action on $\k G\ot\int^l_A$ is diagonal.
\end{enumerate}
\end{prop}

\section{Homological integrals of Calabi-Yau Hopf
algebras}\label{sec2}

In this section we study Noetherian CY Hopf algebras. We show that a
Noetherian CY Hopf algebra has trivial homological integrals, and
its antipode must be bijective.

Let $A$ be an algebra. Recall that $A$ is said to be a {\it
Calabi-Yau algebra of dimension $d$} (cf. \cite{G,BT}) if (i) $A$ is
homologically smooth, that is; $A$ has a bounded resolution of
finitely generated projective $A$-$A$-bimodules, (ii)
$\Ext^i_{A^e}(A,A^e)=0$ if $i\neq d$ and $\Ext_{A^e}^d(A,A^e)\cong
A$ as $A$-$A$-bimodules, where $A^e=A\ot A^{op}$ is the enveloping
algebra of $A$. In what follows, Calabi-Yau is abbreviated to CY for
short.

Let $A$ be an algebra, $\sigma:A\to A$ an algebra morphism, and $M$
a right $A$-module. Denote by $M^\sigma$ the right $A$-module
twisted by the algebra morphism $\sigma$. If $N$ is an
$A$-$A$-bimodule, we denote by ${}^1N^\sigma$ the bimodule whose
right $A$-action is twisted by $\sigma$.

Let $H$ be a Hopf algebra with antipode $S$. We write $\k$ or $\k_H$
for the trivial module defined by the counit of $H$. Let $M$ be an
$H$-$H$-bimodule. Denote $M^{ad}$ the left adjoint $H$-module
defined by $h\cdot m=\sum_{(h)}h_{(1)}mS(h_{(2)})$ for $h\in H$ and
$m\in M$.

Let $D:H\to H\ot H^{op}$ be the map defined by $D=(1\ot
S)\circ\Delta$. Then $D$ is an algebra morphism, and $H\ot H^{op}$
is a free left (and a free right) $H$-module through $D$ (see
\cite[Sect.2]{BZ}). We write $L(H\ot H^{op})$ (resp. $R(H\ot
H^{op})$) for the left (resp. right) $H$-module defined through $D$.
Let ${}_\bullet H\ot H$ be the left $H$-module defined by the left
multiplication of $H$ to the left factor, and $H_\bullet\ot H$ be
the free right $H$-module defined by the right multiplication of $H$
to the left factor. Then $L(H\ot H^{op})\cong {}_\bullet H\ot H$ and
$R(H\ot H^{op})\cong H_\bullet\ot H$. The isomorphisms are given as
follows:
\begin{equation}\label{tag23}
    \varphi:L(H\ot H^{op})\to {}_\bullet H\ot H,\ g\ot h\mapsto
\sum_{(g)}g_{(1)}\ot hS^2(g_{(2)}),
\end{equation}
with its inverse
\begin{equation}\label{tag24}\phi:{}_\bullet H\ot H\to L(H\ot H^{op}),\ g\ot h\mapsto \sum_{(g)}g_{(1)}\ot
hS(g_{(2)});
\end{equation}and
\begin{equation}\label{tag25}
    \psi:R(H\ot H^{op})\to H_\bullet\ot H,\ g\ot
h\mapsto\sum_{(g)}g_{(1)}\ot g_{(2)}h,
\end{equation} with its inverse
\begin{equation}\label{tag26}
\xi:H_\bullet\ot H\to R(H\ot H^{op}),\ g\ot
h\mapsto\sum_{(g)}g_{(1)}\ot S(g_{(2)})h.
\end{equation}

Clearly, $L(H\ot H^{op})$ is an $H$-$H^e$-bimodule and $R(H\ot
H^{op})$ is an $H^e$-$H$-bimodule. Let $M$ be an $H$-$H$-bimodule.
Then one has $M^{ad}\cong\Hom_{H^e}(R(H\ot H^{op}),M)$. By
\cite[Lemma 2.2]{BZ}, the functor $(-)^{ad}$ preserves injective
modules. This property implies the key fact that the Hochschild
cohomology of a bimodule $M$ over a Hopf algebra $H$ can be computed
through the extension groups of the trivial module ${}_H\k$ by
$M^{ad}$, see \cite[Lemma 2.4]{BZ} or \cite[Prop. 5.6]{GK}:

\begin{lem} Let $H$ be a Hopf algebra, and $M$ be an $H$-$H$-bimodule.
Then we have $\Ext_{H^e}^i(H,M)\cong\Ext^i_H({}_H\k,M^{ad})$ for all
$i$.
\end{lem}

Let $H$ be an AS-Gorenstein Hopf algebra. The left homological
integrals $\int^l_H$ of $H$ is a one-dimensional right $H$-module,
and the $H$-module structure is defined through an algebra morphism
$\pi:H\to\k$. We have an algebra automorphism $\nu:H\to H$ defined
by $\nu(h)=\sum_{(h)}\pi(h_{(1)})h_{(2)}$ for $h\in H$. Then, as
right $H$-modules, $\int^l_H\cong \k^\nu$. The following corollary
is proved in \cite[Prop. 4.5]{BZ}. We include the proof for the
completeness here. Notice that the hypothesis that the antipode is
bijective is dropped.

\begin{cor}\label{cor2} Let $H$ be a Noetherian AS-Gorenstein Hopf algebra with injective dimension $d$.
Then $\Ext^i_{H^e}(H,H^e)=0$ for $i\neq d$ and
$\Ext^d_{H^e}(H,H^e)={}^1H^{S^2\nu}$.
\end{cor} \proof The proof is just a slight modification of that of \cite[Prop. 4.5]{BZ}.
Since $H$ is noetherian, we have
$$\begin{array}{ccl}
    \Ext_{H^e}^i(H,H^e)&\cong&\Ext^i_H({}_H\k,L(H\ot H^{op}))\\
    &\cong&\Ext^i_H({}_H\k,H)\ot_HL(H\ot
H^{op}).
  \end{array}
$$ Hence $\Ext^i_{H^e}(H,H^e)=0$ for $i\neq d$, and
$$\Ext^d_{H^e}(H,H^e)\cong\k^\nu\ot_HL(H\ot
H^{op})\overset{(a)}\cong \k^\nu\ot_H({}_\bullet H\ot
H)\overset{(b)}\cong{}^1H^{S^2\nu},$$ where $\nu$ is the algebra
automorphism of $H$ corresponding to the left homological integrals
$\int^l_H$. The isomorphism $(a)$ is given by the isomorphism
$\varphi$ constructed through the map (\ref{tag23}) above; and the
isomorphism $(b)$ holds because the right $H^e$-module structure on
${}_\bullet H\ot H$ induced by the isomorphism $\varphi$ is given as
$(g\ot h)\cdot(x\ot y)=\sum_{(x)}gx_{(1)}\ot yhS^2(x_{(2)})$ for
$g,h,x,y\in H$. \qed

Now we arrive at the main result of this section. Recall from
\cite{LWZ} that an AS-Gorenstein Hopf algebra is {\it unimodular} if
$\int^l_H\cong\k_H$ as right $H$-modules.

\begin{thm}\label{prop4} Let $H$ be a Noetherian Hopf algebra. Then $H$ is CY of dimension $d$ if and only
if

{\rm(i)} $H$ is AS-regular with global dimension $\text{\rm
gldim}(H)=d$ and unimodular,

{\rm(ii)} $S^2$ is an inner automorphism of $H$.
\end{thm} \proof Suppose that $H$ is CY of dimension $d$. By \cite[Lemma
4.1]{K} the triangulated category $D^b_{fd}(H)$ is a CY category of
dimension $d$, where $D^b_{fd}(H)$ is the full triangulated
subcategory of the derived category of $H$ formed by complexes whose
homology is of finite total dimension. Hence gl$\dim H=d$. It is
well-known that, for $i\ge0$, $$\Ext_H^i({}_H\k
,H)\cong\Ext^i_{H^e}(H,\Hom_\k({}_H\k ,H))\cong \Ext^i_{H^e}(H,H\ot
\k_H),$$ where the left $H^e$-bimodule structure on $H\ot \k_H$ is
given by the left multiplication of $H$ on the first factor and the
right $H$-action on the trivial module $\k_H$. Since $H$ is CY of
dimension $d$, we may choose a finitely generated projective
resolution of the $H^e$-module $H$ as follows
$$P^\bullet:\qquad 0\longrightarrow
P^{-d}\longrightarrow\cdots\longrightarrow P^{-1}\longrightarrow
P^0\longrightarrow H\longrightarrow0.$$  Then we have isomorphisms
of complexes
\begin{equation}\label{tag12}\begin{array}{ccl}
\Hom_{H^e}(P^\bullet,H\ot
\k_H)&\cong&\Hom_{H^e}(P^\bullet,H^e)\ot_{H^e} (H\ot\k_H)\\
&\cong&\k\ot_H\Hom_{H^e}(P^\bullet,H^e)\ot_HH\\
&\cong&\k\ot_H\Hom_{H^e}(P^\bullet,H^e).
  \end{array}
\end{equation} Let $Q^\bullet:=\Hom_{H^e}(P^\bullet,H^e)$. Since $H$ is CY of
dimension $d$, $Q^\bullet[d]$ is a projective resolution of the
$H^e$-module $H$. Hence we have
$$\k\ot_H\Hom_{H^e}(P^\bullet,H^e)=\k\ot_HQ^\bullet\overset{\simeq}\longrightarrow\k_H[-d],$$
where the second map is a quasi-isomorphism of right $H$-modules.
Note that the isomorphisms in (\ref{tag12}) are also right
$H$-module morphisms. By taking the cohomology of the complexes in
(\ref{tag12}), we obtain the statement (i).

According to Part (i) and Corollary \ref{cor2},
$\Ext^d_{H^e}(H,H^e)\cong {}^1H^{S^2}$. On the other hand, the CY
property of $H$ implies $\Ext^d_{H^e}(H,H^e)\cong H$ as
$H$-$H$-bimodules. Hence $H$ and ${}^1H^{S^2}$ are isomorphic as
$H$-$H$-bimodules. Therefore $S^2$ must be an inner automorphism.

Conversely, since $H$ is Noetherian and of finite global dimension,
by \cite[Lemma 5.2]{BZ} $H$ is homologically smooth. The assertion
(i) and Corollary \ref{cor2} insure $\Ext^i_{H^e}(H,H^e)=0$ and
$\Ext^d_{H^e}(H,H^e)\cong {}^1H^{S^2}$. The assertion (ii) tells us
that ${}^1H^{S^2}$ is isomorphic to $H$ as an $H$-$H$-bimodule.\qed

\section{Group actions on universal enveloping algebras}\label{sec3}

In this section we consider the universal enveloping Hopf algebra of a Lie algebra and study the CY property of its smash product Hopf algebra.
Let $\mathfrak{g}$ be a finite dimensional Lie algebra,
and $U(\mathfrak{g})$  the universal enveloping algebra of
$\mathfrak{g}$. Recall from \cite[Prop.6.3]{BZ} that $U(\mathfrak{g})$ is an AS-regular Hopf algebra. Now let $G$ be a finite group. We say that
$\mathfrak{g}$ is a left {\it $G$-module Lie algebra} if there is a
$G$-action on $\mathfrak{g}$ such that $\mathfrak{g}$ is a left
$G$-module and $g[x,y]=[gx,gy]$ for all $g\in G$ and $x,y\in
\mathfrak{g}$. If $\mathfrak{g}$ is a left $G$-module Lie algebra,
then $U(\mathfrak{g})$ is an augmented left $G$-module algebra. For
a Lie algebra $\mathfrak{g}$, we write $Aut_{Lie}(\mathfrak{g})$ to
be the group of Lie algebra automorphisms. If $\mathfrak{g}$ is a
left $G$-module Lie algebra, write the associated group morphism as
$\nu:G\longrightarrow Aut_{Lie}(\mathfrak{g})$.


Assume $\dim \mathfrak{g}=d$. Consider the Chevalley-Eilenberg
resolution of the trivial $U(\mathfrak{g})$-module (cf. \cite[Ch.
8]{CE} or \cite[Ch. 10]{L}): {\small\begin{equation}\label{tag13}
0\longrightarrow U(\mathfrak{g})\ot
\wedge^d\mathfrak{g}\overset{\partial^d}\longrightarrow\cdots\overset{\partial^3}\longrightarrow
U(\mathfrak{g})\ot \mathfrak{g}\wedge
\mathfrak{g}\overset{\partial^2}\longrightarrow U(\mathfrak{g})\ot
\mathfrak{g}\overset{\partial^1}\longrightarrow
U(\mathfrak{g})\longrightarrow {}_{U(\mathfrak{g})}\k
\longrightarrow0,
\end{equation}}
 where for $x_1,\dots,x_n\in \mathfrak{g}$, $$\begin{array}{cl}
            & \partial^n(1\ot x_1\wedge\cdots\wedge x_n)\\
             = &\displaystyle \sum_{i=1}^n(-1)^{i+1}x_i\ot x_1\wedge\cdots\wedge
\hat{x}_i\wedge\cdots\wedge x_n \\
               & +\displaystyle\sum_{1\leq i<j\leq n}(-1)^{i+j}1\ot
              [x_i,x_j]\wedge x_1\wedge\cdots\wedge
              \hat{x}_i\wedge\cdots\wedge \hat{x}_j\wedge\cdots\wedge
              x_n.
           \end{array}
$$
Since $\mathfrak{g}$ is a left $G$-module, $\wedge^n\mathfrak{g}$ is
a left $G$-module with the diagonal action. Thus $U(\mathfrak{g})\ot
\wedge^n\mathfrak{g}$ is a left
$U(\mathfrak{g})\#\k G$-module. It is not hard to check that
the differentials in the resolution above are also left $G$-module
maps. Hence the resolution above is in fact a projective resolution
of the left $U(\mathfrak{g})\#\k G$-module
${}_{U(\mathfrak{g})\#\k G}\k $.

\begin{lem}\label{lem1} Let $G$ be a finite group, and $\mathfrak{g}$ be a $G$-module Lie algebra of
dimension $d$. Then $U(\mathfrak{g})$ is AS-regular of global
dimension $d$ and as left $G$-modules $\int^l_{U(\mathfrak{g})}\cong
\wedge^d \mathfrak{g}^*$, where left $G$-module action on
$\mathfrak{g}^*$ is defined by $(g\cdot \beta)(x)=\beta(g^{-1}x)$
for $g\in G$, $\beta\in \mathfrak{g}^*$ and $x\in \mathfrak{g}$, and
$G$ acts on $\wedge^d \mathfrak{g}^*$ diagonally.
\end{lem}
\proof Since $\mathfrak{g}$ is of dimension $d$, the universal
enveloping algebra $U(\mathfrak{g})$ has global dimension $d$ (cf.
\cite[Ch. VIII]{CE}). The AS-regularity of $U(\mathfrak{g})$ is
proved in \cite[Prop. 6.3]{BZ}. Applying the functor
$\Hom_{U(\mathfrak{g})}(-,U(\mathfrak{g}))$ to the projective
resolution (\ref{tag13}) of ${}_{U(\mathfrak{g})}\k $ above, we
obtain that $\int^l_{U(\mathfrak{g})}$ is the homology at the final
position of the following complex of left $G$- and right
$U(\mathfrak{g})$-modules (warning: they are not
$G$-$U(\mathfrak{g})$-bimodules)
$$\begin{array}{c}
    0\longrightarrow
\Hom_{U(\mathfrak{g})}(U(\mathfrak{g}),U(\mathfrak{g}))\overset{{\partial^1}^*}\longrightarrow\Hom_{U(\mathfrak{g})}(U(\mathfrak{g})\ot
\g,U(\mathfrak{g}))\overset{{\partial^2}^*}\longrightarrow\cdots  \\
\overset{{\partial^d}^*}
\longrightarrow\Hom_{U(\mathfrak{g})}(U(\mathfrak{g})\ot
\wedge^d\mathfrak{g},U(\mathfrak{g}))\longrightarrow0,
  \end{array}
$$ which is isomorphic to the
following complex of left $G$- and right $U(\mathfrak{g})$-modules
(also not bimodules)
\begin{equation}\label{tag8}
    0\longrightarrow U(\mathfrak{g})\overset{\delta^0}\longrightarrow \g^*\ot U(\mathfrak{g})
\overset{\delta^1}\longrightarrow \cdots\overset{\delta^{d-1}}
\longrightarrow \wedge^d\mathfrak{g}^*\ot
U(\mathfrak{g})\longrightarrow 0.
\end{equation}
Note that the left $G$-action on $\wedge^i\mathfrak{g}^*\ot U(\mathfrak{g})$ is diagonal and $\wedge^i\mathfrak{g}^*\ot U(\mathfrak{g})$ as a right $U(\mathfrak{g})$-module is free. The differential $\delta^{i-1}$ is induced by ${\partial^i}^*$ (for $i\ge1$) through the obvious isomorphisms
$\Hom_{U(\mathfrak{g})}(U(\mathfrak{g})\ot
\wedge^i\mathfrak{g},U(\mathfrak{g}))\cong \wedge^i\mathfrak{g}^*
\ot U(\mathfrak{g})$.

Now let $\{x_1,\dots,x_d\}$ be a basis of
$\mathfrak{g}$ and $\{x_1^*,\dots,x_d^*\}$ the dual basis of
$\mathfrak{g}^*$. Note that the differentials in the complex
(\ref{tag8}) are also right $U(\mathfrak{g})$-module morphisms. The
image of the element $\alpha=x_1^*\wedge\cdots\wedge x_d^*\ot 1$ in
the $d$th cohomology is nonzero. Otherwise it would imply that the
$d$th cohomology is zero. For $\beta\in \wedge^d\mathfrak{g}^*\ot
U(\mathfrak{g})$, let $\overline{\beta}$ be the image of $\beta$ in
the $d$th cohomology. Now for $g\in G$,
$g\cdot\overline{\alpha}=\overline{g\cdot\alpha}=\overline{g\cdot(x_1^*\wedge\cdots\wedge
x_d^*)\ot}1=\omega\overline{\alpha}$, for some nonzero element
$\omega\in k$. Thus we obtain an isomorphism of left
$G$-modules: $\int^l_{U(\mathfrak{g})}\overset{\cong}\longrightarrow
\wedge^d\mathfrak{g}^*$ sending $\overline{\alpha}$ to
$x_1^*\wedge\dots\wedge x_d^*$.  \qed

Let $G$ be a group, $\mathfrak{g}$ a left $G$-module Lie algebra.
It is well known that $U(\mathfrak{g})\#\k G$ is a
cocommutative Hopf algebra with the coproduct and counit given by
those of $U(\mathfrak{g})$ and of $\k G$.

\begin{lem}\label{lem3} Let $G$ be a finite group. If $\mathfrak{g}$ is a finite dimensional $G$-module Lie algebra,
then $U(\mathfrak{g})\#\k G$ is an AS-regular algebra.
\end{lem}
\proof Since $\mathfrak{g}$ is finite dimensional, $U(\mathfrak{g})$
is an AS-regular algebra. Hence $U(\mathfrak{g})\#\k G$ has
finite global dimension. Now the statement follows from Propositions
\ref{prop1} and \ref{prop2}. \qed

\begin{lem}\label{lem4} Let $G$ be a finite group, and $\mathfrak{g}$ a finite dimensional left
$G$-module Lie algebra. If $U(\mathfrak{g})\#\k G$ is CY of
dimension $d$, then $U(\mathfrak{g})$ is CY of dimension $d$.
\end{lem}
\proof Write $B$ for $U(\mathfrak{g})\#\k G$. By Lemma \ref{lem3},
$B$ is AS-regular of global dimension $d$. By Proposition
\ref{prop1}, $\int^l_B\cong(\int^l_{U(\mathfrak{g})}\#\k G)^G$.
Choose a basis $\alpha\#t$ of $\int^l_B$  with
$\alpha\in\int^l_{U(\mathfrak{g})}$ and $t\in \k G$. Since $B$ is a
Noetherian cocommutative Hopf algebra, by Theorem \ref{prop4} the
right $B$-module action on $\int^l_B$ is trivial. It follows that
$\alpha\#t=(\alpha\#t)\cdot(1\#g)=\alpha\#tg$ for all $g\in G$. This
implies $tg=t$ for all $g\in G$ and hence $t$ is an integral of $\k
G$. We may now assume $t=\sum_{g\in G}g$. For $a\in
U(\mathfrak{g})$, we have
$$\varepsilon(a)\alpha\#t=(\alpha\#t)\cdot(a\#1)=\sum_{(t)}\alpha\cdot(t_{(1)}a)\#t_{(2)}=
\sum_{g\in G}\alpha\cdot(ga)\# g,$$
 which forces
$\alpha\cdot(ga)=\varepsilon(a)\alpha$ for all $g\in G$. Replacing
$a$ by $g^{-1}a$, we obtain $\alpha\cdot
a=\varepsilon(g^{-1}a)\alpha=\varepsilon(a)\alpha$ for all $a\in
U(\mathfrak{g})$. Therefore, the right $U(\mathfrak{g})$-action on
the integral space $\int_{U(\mathfrak{g})}^l$ is trivial. Now the
result follows from Theorem \ref{prop4}. \qed

\begin{thm}\label{thm1} Let $\mathfrak{g}$ be a finite dimensional Lie
algebra, and $G\subseteq Aut_{Lie}(\mathfrak{g})$  a finite group.
Then the skew group algebra $U(\mathfrak{g})\#\k G$ is a CY algebra
of dimension $d$ if and only if $U(\mathfrak{g})$ is a CY algebra of
dimension $d$ and $G\subseteq SL(\mathfrak{g})$.
\end{thm}
\proof Suppose $G\subseteq SL(\mathfrak{g})$ and $U(\mathfrak{g})$
is CY. As before, write $B$ for $U(\mathfrak{g})\#\k G$. Since $B$
is a cocommutative Hopf algebra, by Theorem \ref{prop4} we only need
to show that $\int^l_B\cong \k _B$ as right $B$-modules, where $\k
_B$ is the trivial right $B$-module. From Proposition \ref{prop1},
we have $\int^L_B\cong(\int^l_{U(\mathfrak{g})}\#\k G)^G$. By Lemma
\ref{lem1}, $\int^l_{U(\mathfrak{g})}\cong\wedge^d\mathfrak{g}^*$ as
left $G$-modules. Let $\{x_1^*,\dots,x_d^*\}$ be a basis of $\g^*$.
We have $g\cdot(x_1^*\wedge\dots\wedge
x_d^*)=\det(g)^{-1}x_1^*\wedge\dots\wedge
x_d^*=x_1^*\wedge\dots\wedge x_d^*$ for $g\in G\subseteq
SL(\mathfrak{g})$. If $\alpha$ is a basis of
$\int^l_{U(\mathfrak{g})}$, then $g\cdot \alpha=\alpha$ for all
$g\in G$. Assume that $\alpha\# t$ is an element in
$(\int^l_{U(\mathfrak{g})}\#\k G)^G$ for some $t\in \k G$. Then
$\alpha\# t=g\cdot(\alpha\# t)=g\cdot\alpha\# gt=\alpha\#gt$ for all
$g\in G$. So we have $gt=t$ for all $g\in G$. Hence $t$ must be an
integral of $\k G$. Now we may assume that $\alpha\# t$ is a basis
of $\int^l_B$ with $t$ a nonzero integral of $\k G$. By assumption,
$U(\mathfrak{g})$ is a cocommutative CY Hopf algebra. It follows
from Theorem \ref{prop4} that the right $U(\mathfrak{g})$-module
structure on $\int^l_{U(\mathfrak{g})}$ is trivial. Now for $a\in
U(\mathfrak{g})$ and $g\in G$, $(\alpha\ot t)\cdot
(a\#g)=\sum_{(t)}\alpha\cdot (t_{(1)}a)\ot
t_{(2)}g=\sum_{(t)}\alpha\ot\varepsilon(t_{(1)}a)t_{(2)}g=\alpha\ot\varepsilon(a)tg=\varepsilon(a)\alpha\ot
t=\varepsilon_B(a\#g)\alpha\ot t$. Thus the right $B$-module
structure of $\int^l_B$ is trivial. Therefore $B$ is a CY Hopf
algebra of dimension $d$.

Conversely, assume $B$ is a CY algebra. By Lemma \ref{lem4},
$U(\mathfrak{g})$ is CY. By the proof of Lemma \ref{lem4}, we may
assume that $\alpha\ot t$ is a basis of $\int^l_B$ with $t$ an
integral of $\k G$ and $\alpha\in\int^l_{U(\g)}$. Note that
$\int^l_B\cong(\int^l_{U(\mathfrak{g})}\# \k G)^G$. Hence for $g\in
G$, $\alpha\ot t=g\cdot (\alpha\ot t)=g\cdot\alpha\ot
gt=g\cdot\alpha\ot t$. We get $g\cdot\alpha=\alpha$. This implies
that
 the left $G$-action on $\int^l_{U(\mathfrak{g})}=\wedge^d\mathfrak{g}^*$
is trivial. Thus  we  have $\det(g)=1$ for all $g\in G$, i.e.,
$G\subseteq SL(\mathfrak{g})$. \qed

Now let $G$ be an arbitrary finite group, $\mathfrak{g}$ a finite
dimensional $G$-module Lie algebra. As before, we let
$\nu:G\longrightarrow Aut_{Lie}(\mathfrak{g})$ be the associated
group map. From the proof of the theorem above, we obtain

\begin{cor}\label{cor1} Let $G$ and $\mathfrak{g}$ be as above. Then $U(\mathfrak{g})\#\k G$
is CY if and only if $U(\mathfrak{g})$ is CY and $\im(\nu)\subseteq
SL(\mathfrak{g})$.
\end{cor}

\section{Cocommutative CY Hopf algebras of low dimensions}

Let $A$ be an augmented algebra with a fixed augmentation map
$\varepsilon:A\longrightarrow\k $. If $A$ is a CY algebra of
dimension $d$, then by \cite[Lemma 4.1]{K} the shift functor $[d]$
of the triangulated category $D^b_{fd}(A)$ is a (graded) Serre
functor (see the appendix of \cite{B}), where $D^b_{fd}(A)$ is the
full triangulated subcategory of the derived category of $A$
consisting of complexes with finite dimensional total cohomologies.
Hence
$\Ext_A^d({}_A\k ,{}_A\k )\cong\Ext_A^0({}_A\k ,{}_A\k )\cong\k $.

Now let $A=U(\mathfrak{g})$ be the universal enveloping algebra of a
Lie algebra $\mathfrak{g}$ of dimension $d$. Tensoring with
$\k _A$ the Chevalley-Eilenberg resolution
of ${}_A\k $, we obtain the following complex:
\begin{equation}\label{tag10}0\longrightarrow
\wedge^d\mathfrak{g}\overset{\delta^d}\longrightarrow\cdots\overset{\delta^3}\longrightarrow
\mathfrak{g}\wedge \mathfrak{g}\overset{\delta^2}\longrightarrow
\mathfrak{g}\overset{\delta^1}\longrightarrow
\k \longrightarrow0,
\end{equation} where the differential is given as, for $2\leq n\leq d$, and
$x_1,\dots,x_n\in \mathfrak{g}$,
 $$   \begin{array}{cl}
            & \delta^n(x_1\wedge\cdots\wedge x_n)\\
             &=\displaystyle\sum_{1\leq i<j\leq n}(-1)^{i+j}
              [x_i,x_j]\wedge x_1\wedge\cdots\wedge
              \hat{x}_i\wedge\cdots\wedge \hat{x}_j\wedge\cdots\wedge
              x_n,
           \end{array}
$$
and $\delta^1=0$. The $n$th homology of the
           complex above is
           $\Tor^A_n(\k _A,{}_A\k )$.

The following lemma  can be deduced from \cite[Proposition 6.3]{BZ}.

\begin{lem}\label{lem5} Let $\mathfrak{g}$ be a Lie algebra of dimension $d$.
The following are equivalent.

\begin{enumerate}
\item[(i)] $U(\mathfrak{g})$ is a CY algebra.
\item[(ii)] $\Ext_{A}^d({}_{A}\k_{A}\k )\neq0$.
\item[(iii)] The differential $\delta^d$ in the complex (\ref{tag10})
is zero.
\item[(iv)] $\text{\rm tr}(\text{\rm
ad}_\g(x))=0$ for all $x\in\g$.
\end{enumerate}
\end{lem}
\proof (i) $\Rightarrow$ (ii) $\Rightarrow $ (iii) are obvious. (iv) $\Rightarrow$ (i) follows from \cite[Proposition 6.3]{BZ}. We just need to show that (iii) $\Rightarrow $ (iv). Assume $\delta^d=0$. Let
$\{x_1,\dots,x_d\}$ be a basis of $\mathfrak{g}$. Then
$$
\begin{array}{ccl}
 0 & = &\delta^d(x_1\wedge\cdots\wedge x_d)  \\
 & = &\displaystyle \sum_{1\leq i<j\leq d}(-1)^{i+j}[x_i,x_j]\wedge
x_1\wedge\cdots\wedge\hat{x}_i\wedge\cdots\wedge\hat{x}_j\wedge\cdots\wedge
x_d \\
&=&\displaystyle\sum_{i=1}^d(-1)^i\text{tr}(\text{ad}_\mathfrak{g}(x_i))x_1\wedge\cdots\wedge\hat{x}_i\wedge\cdots\wedge x_d.
\end{array}
$$
 We have $\text{tr}(\text{ad}_\mathfrak{g}(x_i))=0$ for all
$1\leq i\leq d$. \qed


Let $H$ be a Hopf algebra. Denote by $P(H)$ the space of all
 primitive elements of $H$ and denote by $G(H)$  the group of all
group-like elements of $H$.

\begin{thm}\label{thm2} Let $H$ be a cocommutative Hopf algebra such that
$\dim P(H)<\infty$ and $G(H)$ is finite. Then $H$ is CY of dimension
2 if and only if there is a finite group $G$ and a group map
$\nu:G\to SL(2,\k )$ such that $H\cong\k [x,y]\#\k G$, where the $G$-action on
$\k [x,y]$ is given by $\nu$.
\end{thm}
\proof The sufficiency follows from Corollary \ref{cor1}. For the
necessity, it is well-known that a cocommutative Hopf algebra $H$
over an algebraic closed field is isomorphic to $U(\mathfrak{g})\#\k
G$, where $\mathfrak{g}=P(H)$ is the Lie algebra of all the
primitive elements of $H$ and $G=G(H)$ is the group of the
group-like elements of $H$ (cf. \cite{La,Sw}). Since $H$ is CY of
dimension $2$, by Lemma \ref{lem4}, $U(\mathfrak{g})$ is CY of
dimension 2. Hence we get the global dimension of $U(\mathfrak{g})$
is 2, which implies $\dim(\mathfrak{g})=2$. By Lemma \ref{lem5},
$\mathfrak{g}$ must be abelian. Hence $U(\mathfrak{g})\cong\k
[x,y]$. The rest of the proof follows from Corollary \ref{cor1}.\qed

\begin{rem} {\rm Let $K=\im(\nu)\subseteq SL(2,\k )$ and
$N=\ker(\nu)$. Then $N$ is a normal subgroup of $G$. There is a weak
$H$-action on the group algebra $\k N$ (cf. \cite{BCM}), and
there is a map $\sigma:K\times K\longrightarrow \k N$ so
that $\k G$ is isomorphic to the crossed product of
$\k N$ and $K$, that is; $\k G\cong
\k N\#_\sigma \k K$ (cf. \cite{BCM} or \cite[Ch.
4]{P}). Since $K\subseteq SL(2,\k )$, $K$ acts naturally on
$\k [x,y]$. So $K$ acts on $\k [x,y]\ot
\k N$ diagonally. The map $\sigma:K\times K\longrightarrow
\k N$ may be extended to $K\times K\longrightarrow
\k [x,y]\ot \k N$ (the map is also denoted by
$\sigma$). Then one can check a CY cocommutative pointed Hopf
algebra as in Theorem \ref{thm2} has the form  $H\cong
(\k [x,y]\ot\k N)\#_\sigma\k K.$ }\end{rem}

Next we discuss 3-dimensional CY cocommutative Hopf algebras.
We know that a cocommutative Hopf algebra is the skew-group algebra
of a universal enveloping algebra of a Lie algebra and a group
algebra. Let us now deal with 3-dimensional Lie algebras. By
Lemma \ref{lem5}, we may list all 3-dimensional Lie algebras
whose universal enveloping algebra is CY, since the 3-dimensional
Lie algebras are classified completely. However, let us first get a
view of the Lie bracket of such Lie algebras over an arbitrary
basis. Let $\mathfrak{g}$ be a 3-dimensional vector space with a
basis $\{x,y,z\}$. Define a bracket on $\mathfrak{g}$ as follows:
\begin{eqnarray}
 \nonumber [x,y] & = & ax+by+wz, \\
  \label{tag11}[x,z] & = & cx+vy-bz, \\
 \nonumber [y,z] & = & ux-cy+az,
\end{eqnarray}
where $a,b,c,u,v,w\in\k $. A direct verification shows that $\mathfrak{g}$ is a
Lie algebra.

\begin{lem}\label{lem6} With the bracket defined above, $\mathfrak{g}$ is a
Lie algebra.
\end{lem}

Now we have the following easy but useful result.

\begin{prop}\label{prop3} Let $\mathfrak{g}$ be a 3-dimensional Lie
algebra, and $\{x,y,z\}$ be a basis of $\g$. Then $U(\mathfrak{g})$
is a CY algebra if and only if the Lie bracket is given by
{\rm(\ref{tag11})}.
\end{prop}
\proof The lemma above shows that the bracket defined in
(\ref{tag11}) is a Lie bracket on $\mathfrak{g}$. By Lemma
\ref{lem5}, $U(\mathfrak{g})$ is a CY algebra if
tr$(\text{ad}_\mathfrak{g}(l))=0$ for all $l\in \mathfrak{g}$. Now
it is direct to check tr$(\text{ad}_\mathfrak{g}(x))=0$,
tr$(\text{ad}_\mathfrak{g}(y))=0$ and
tr$(\text{ad}_\mathfrak{g}(z))=0$. It follows that $U(\mathfrak{g})$
is a CY algebra.

Suppose $U(\mathfrak{g})$ is a CY algebra. Assume
$$\begin{array}{ccl}
[x,y] & = & k_{11}x+k_{12}y+k_{13}z, \\
{[x,z]} & = & k_{21}x+k_{22}y+k_{23}z, \\
{[y,z]} & = & k_{31}x+k_{32}y+k_{33}z,
         \end{array}
$$ where $k_{ij}\in\k $. By Lemma \ref{lem5}, the
differential $\delta^d$ in the complex (\ref{tag10}) associated to
$\mathfrak{g}$ is zero, that is; $\delta^3(x\wedge y\wedge z)=0$.
Then we get $[x,y]\wedge z-[x,z]\wedge y+[y,z]\wedge x=0$, which is
equivalent to $$(k_{33}-k_{11})x\wedge z-(k_{12}-k_{23})y\wedge
z+(k_{21}+k_{32})x\wedge y=0.$$ Hence we get $k_{11}=k_{33}$,
$k_{12}=-k_{23}$ and $k_{21}=-k_{32}$. Therefore the Lie bracket is
of the form given in (\ref{tag11}).\qed

Let $\mathfrak{g}$ be a 3-dimensional vector space. Fix a basis
$\{x,y,z\}$ of $\mathfrak{g}$. Proposition \ref{prop3} states that,
given a sextuple $(a,b,c,u,v,w)\in\k ^6$, there is a Lie bracket on
$\mathfrak{g}$ defined by this sextuple via (\ref{tag11}) so that
the universal enveloping algebra of $\mathfrak{g}$ is a
3-dimensional CY algebra. Moreover, any 3-dimensional CY universal
enveloping algebra is obtained in this way. We list below all
3-dimensional Lie algebras whose universal enveloping algebras are
CY.

\begin{prop}\label{prop5} Let $\mathfrak{g}$ be a finite dimensional Lie algebra.
Then $U(\mathfrak{g})$ is CY of dimension 3 if and only if
$\mathfrak{g}$ is isomorphic to one of the following Lie algebras:
\begin{itemize}
  \item [(i)] The 3-dimensional simple Lie algebra $\mathfrak{sl}(2,\k )$;
  \item [(ii)] $\mathfrak{g}$ has a basis $\{x,y,z\}$ such that $[x,y]=y$,
$[x,z]=-z$ and $[y,z]=0$.
  \item [(iii)] The Heisenberg algebra, that is; $\mathfrak{g}$ has a
basis $\{x,y,z\}$ such that $[x,y]=z$ and $[x,z]=[y,z]=0$;
\item[(iv)] The 3-dimensional abelian Lie algebra.
\end{itemize}
\end{prop}\proof  Note that by Proposition \ref{prop3} the universal enveloping algebras of the Lie algebras listed above
are CY of dimension 3. We show that they are the only
possible cases. We divide the 3-dimensional Lie algebras
into 4 classes: \\
Case 1: dim$[\mathfrak{g},\mathfrak{g}]=3$, that is, $\mathfrak{g}=[\mathfrak{g},\mathfrak{g}]$;\\
Case 2: dim$[\mathfrak{g},\mathfrak{g}]=2$;\\
Case 3: dim$[\mathfrak{g},\mathfrak{g}]=1$;\\
Case 4: dim$[\mathfrak{g},\mathfrak{g}]=0$ or $\mathfrak{g}$ is abelian.

Case 1. If $\mathfrak{g}=[\mathfrak{g},\mathfrak{g}]$, then it is
well-known that $\mathfrak{g}\cong\mathfrak{sl}(2,\mathbb{C})$. This
gives us the Lie algebra (i).

Case 2. Assume that the Lie algebra $\mathfrak{g}$ has
$\dim[\mathfrak{g},\mathfrak{g}]=2$. We choose a proper
basis $\{x,y,z\}$ for $\mathfrak{g}$ so that $\mathfrak{g}$ satisifies (cf. \cite{EW}):\\

(a) $[x,y]=y$, $[x,z]=\mu z$ and $[y,z]=0$, where
$0\neq\mu\in\mathbb{C}$; or

(b) $[x,y]=y$, $[x,z]=y+ z$ and $[y,z]=0$.

Since $\mathfrak{g}$ is CY, it follows from Proposition \ref{prop3} that the Lie bracket of $\mathfrak{g}$ must satisfy the relations in (\ref{tag11}). In the case (a), we must have $\mu=-1$. So $\mathfrak{g}$ is the Lie algebra
given by (ii).  Since the defining relations in the case (b) do not satisfy
(\ref{tag11}), the Lie bracket defined in (b) does not define a Lie
algebra with CY universal enveloping algebra.

Case 3. Assume that the Lie algebra $\mathfrak{g}$ has
$\dim[\mathfrak{g},\mathfrak{g}]=1$. Similar to Case 2, by choosing a proper
basis $\{x,y,z\}$, we see that $\mathfrak{g}$ is determined by either of the following two cases:

(a) $[\mathfrak{g},\mathfrak{g}]$ is contained in the center of
$\mathfrak{g}$. In this case, $\mathfrak{g}$ is the Heisenberg algebra: $[x,y]=z$, and $[x,z]=[y,z]=0$.

(b) $[\mathfrak{g},\mathfrak{g}]$ is  not contained in the center of
$\mathfrak{g}$. In this case, we can have: $[x,y]=y$ and
$[x,z]=[y,z]=0$.

Clearly, the Lie algebra defined by the case (b) does not satisfy
the relations in (\ref{tag11}), and hence its universal enveloping
algebra can not be CY. Therefore, we have only the Heisenberg Lie algebra (iii).

Case 4. When a 3-dimensional Lie algebra is abelian, then its
universal enveloping algebra is certainly CY. This yields the class
(iv). \qed

Now, similar to Theorem \ref{thm2}, we may write down all possible
3-dimensional Noetherian CY cocommutative Hopf algebras with a
finite number of group-like elements.

\begin{thm}\label{thm5} Let $H$ be a cocommutative Hopf algebra such that $\dim P(H)<\infty$ and $G(H)$ is
finite.  Then $H$ is CY of dimension 3 if and only if $H\cong
U(\g)\# \k G$, where $\g$ is one of the 3-dimensional algebras
listed in Proposition \ref{prop5} and $G$ is a finite group with a
group morphism $\nu:G\to Aut_{Lie}(\g)$ such that $\im(\nu)$ is also
a subgroup of $SL(\g)$.
\end{thm}
\proof The proof is similar to that of Theorem \ref{thm2}. \qed

\begin{rem} {\rm The cocommutative Hopf algebra discussed in
this section possesses  finite number of group-like elements. If the group of
group-like elements is infinite, then the situation becomes very complicated. For an infinite group, it is hard to determine when the group algebra is CY,  even in the low dimensional cases.
However, there are some examples of CY group algebras of low
dimensions (see \cite{LWZ}). If $G$ is a finitely
generated group such that $\k G$ is Noetherian and is of
GK-dimension 1, then $G\cong\mathbb Z$ (cf. \cite[Prop. 8.2]{LWZ}). In this case, $\k G$ is CY of dimension 1. Example 8.5 of \cite{LWZ}
provides us an example of Noetherian affine CY group algebra of
dimension $2$.}
\end{rem}

\section{Sridharan enveloping algebras}

In this section, we discuss the CY property of a Sridharan enveloping algebra of a finite dimensional Lie algebra.
In general, a Sridharan enveloping algebra is no longer a  Hopf algebra,
but a cocycle deformation of a cocommutative Hopf algebra or a Poincar\'{e}-Birkhoff-Witt (PBW) deformation of a polynomial algebra (cf. \cite{Sr,N0,N}).
We will see that the CY property of a Sridharan enveloping algebra is closely related to the CY property of a universal enveloping algebra.
The class of Sridharan algebras contains many interesting algebras, such as Weyl algebras. Many homological properties of
Sridharan enveloping algebras have been discussed in
\cite{Sr,Ka,N0,N}, especially the Hochscheld (co)homology and the cyclic
homology. In \cite{N0}, Nuss listed all nonisomorphic Sridharan
enveloping algebras of 3-dimensional Lie algebras. Based on these
results, we obtain in this section some equivalent conditions for a
Sridharan enveloping algebra of a finite dimensional Lie algebra to
be CY. We then list all possible nonisomorphic 3-dimensional CY
Sridharan enveloping algebras, and partly answer a question
proposed by Berger at the end of \cite{Be}.

Let $\g$ be a finite dimensional Lie algebra, and let $f\in
Z^2(\g,\k)$ be an arbitrary 2-cocycle, that is; $f:\g\times\g\to\k$
such that $$f(x,x)=0\ \text{and}\
f(x,[y,z])+f(y,[z,x])+f(z,[x,y])=0$$ for all $x,y,z\in\g$. The
Sridharan enveloping algebra of $\g$ is defined to be the associative
algebra $U_f(\g)=T(\g)/I$, where $I$ is the two-side ideal of
$T(\g)$ generated by the elements
$$x\ot y-y\ot x-[x,y]-f(x,y),\qquad\text{for all $x,y\in \g$}.$$ For
$x\in\g$, we still denote by $x$ its image in $U_f(\g)$. Clearly,
$U_f(\g)$ is a filtered algebra with the associated graded algebra
$gr(U_f(\g))$ being a polynomial algebra. By \cite[Cor. 3.3]{Sr},
there is one to one correspondence between the group of the algebra
automorphisms $\theta:U_f(\g)\to U_f(\g)$ such that the associated
graded map $gr(\theta)$ is the identity,  and the group $Z^1(\g,\k)$
of the first cocycles. Thus given a 1-cocycle $h\in Z^1(\g,\k)$,
there is an algebra automorphism $\xi_f:U_f(\g)\longrightarrow
U_f(\g)$, for any 2-cocycle $f\in Z^2(\g, \k)$, defined by
\begin{equation}\label{tag21}
    \xi_f(x)=x+h(x)
\end{equation} for all $x\in\g$.
When $f=0$, the map (\ref{tag21}) defines an algebra automorphism
$\xi:U(\g)\longrightarrow U(\g)$. Clearly $\xi_f$ has an inverse
given by $\xi_f^{-1}(x)=x-h(x)$ for all $x\in \g$.

In the sequel, we fix a 2-cocycle $f\in Z^2(\g,\k)$ and let  $A=U_f(\g)$, and
$A^e=A\ot A^{op}$. Define a linear map $D:\g\longrightarrow A^e$ by
$D(x)=x\ot 1-1\ot x$ for all $x\in\g$. By \cite{Sr}, $D$ induces an
algebra morphism $U(\g)\longrightarrow A^e$, still denoted by $D$.
Thus $A^e$ can be viewed both as a left and as a right $U(\g)$-module.
Now let
$\xi:U(\g)\to U(\g)$ and $\xi_f:A\to A$ be defined by (\ref{tag21}). Since
${}^{\xi_f}\!\!A^1\ot A^{op}$ is a left $A^e$-module, it is also a left
$U(\g)$-module. We have the following isomorphisms.

\begin{lem}\label{lem21} As left $U(\g)$- and right $A^e$-bimodules,
$${}^\xi\!(A\ot A^{op})\cong {}^1A^{\xi_f^{-1}}\ot A^{op},$$
$${}^{\xi^{-1}}\!\!(A\ot A^{op})\cong {}^1A^{\xi_f}\ot A^{op}.$$
\end{lem}
\proof It is easy to check that the following
diagram of algebra morphisms is commutative.
$$\xymatrix{
  U(\g) \ar[d]_{\xi} \ar[r]^{D} & A\ot A^{op} \ar[d]^{\xi_f\ot id} \\
  U(\g)\ar[r]^{D} & A\ot A^{op}.   }
$$
It follows that ${}^\xi(A\ot A^{op})\cong {}^{\xi_f}\!\!A^1\ot A^{op}$ as left $U(\g)$- and right $A^e$-bimodules.
On the other hand, we have the $U(\g)$-$A^e$-bimodule isomorphism $\xi_f^{-1}\ot id:{}^{\xi_f}\!\!A^1\ot
A^{op}\longrightarrow {}^1\!A^{\xi_f^{-1}}\ot A^{op}$. The composite of the aforementioned two isomorphisms gives us the desired isomorphism.

The second isomorphism in the lemma can be proved similarly. \qed

As $U(\g)$ is a Hopf algebra, the space $\k$ is a trivial
$U(\g)$-$U(\g)$-bimodule. Thus $\k^\xi$ is a right $U(\g)$-module
twisted by the automorphism $\xi$.

\begin{lem}\label{lem22} As right $A^e$-modules, $\k^\xi\ot_{U(\g)}A^e\cong {}^1\!A^{\xi_f}$.
\end{lem}
\proof Consider the exact sequence of right $U(\g)$-modules:
$$0\longrightarrow I\longrightarrow U(\g)\overset{\varepsilon}\longrightarrow
\k\longrightarrow0.$$ Applying the functor
$-\ot_{U(\g)}{}^1U(\g)^\xi$, we obtain the following exact sequence
of right $U(\g)$-modules:
$$0\longrightarrow I^\xi\longrightarrow
U(\g)^\xi\longrightarrow \k^\xi\longrightarrow0.$$ By \cite[Prop.
5.2]{Sr}, $A^e$ is a free $U(\g)$-module on both sides. Tensoring
the above exact sequence with $A^e$, we obtain the following exact
sequence of right $A^e$-modules:
$$0\longrightarrow I^\xi\ot_{U(\g)}A^e\longrightarrow
U(\g)^\xi\ot_{U(\g)}A^e\longrightarrow
\k^\xi\ot_{U(\g)}A^e\longrightarrow0,$$
which is isomorphic to the
following sequence of right $A^e$-modules:
$$0\longrightarrow
I\ot_{U(\g)}{}^{\xi^{-1}}\!\! (A^e)\longrightarrow
U(\g)\ot_{U(\g)}{}^{\xi^{-1}}\!\!(A^e)\longrightarrow
\k\ot_{U(\g)}{}^{\xi^{-1}}\!\!(A^e)\longrightarrow0.$$
By Lemma \ref{lem21}, the sequence above is isomorphic to the following exact
sequence of right $A^e$-modules:
$$0\longrightarrow I\ot_{U(\g)}({}^1\!A^{\xi_f}\ot A^{op})\longrightarrow
U(\g)\ot_{U(\g)}({}^1\!A^{\xi_f}\ot A^{op})\longrightarrow
\k\ot_{U(\g)}({}^1\!A^{\xi_f}\ot A^{op})\longrightarrow0.$$
Hence we obtain the following right $A^e$-module isomorphisms.
$$\k^\xi\ot_{U(\g)}A^e\cong\k\ot_{U(\g)}({}^1\!A^{\xi_f}\ot
A^{op})\cong\frac{{}^1\!A^{\xi_f}\ot A^{op}}{D(I)({}^1\!A^{\xi_f}\ot
A^{op})}.$$
On the other hand, by a right version of the proof of
\cite[Prop. 5.3]{Sr} we have the following exact sequence of right
$A^e$-modules:
$$0\longrightarrow I\ot_{U(\g)}(A\ot A^{op})\longrightarrow A\ot
A^{op}\longrightarrow A\longrightarrow0.$$
Tensoring it with ${}^1\!A^{\xi_f}$ over $A$, we obtain an exact sequence of right $A^e$-modules:
$$0\longrightarrow I\ot_{U(\g)}(A\ot A^{op})\ot_A{}^1\!A^{\xi_f}\longrightarrow (A\ot
A^{op})\ot_A{}^1\!A^{\xi_f}\longrightarrow
A\ot_A{}^1\!A^{\xi_f}\longrightarrow0,$$ which is isomorphic to
$$0\longrightarrow I\ot_{U(\g)}({}^1\!A^{\xi_f}\ot
A^{op})\longrightarrow {}^1\!A^{\xi_f}\ot A^{op}\longrightarrow
{}^1\!A^{\xi_f}\longrightarrow0.$$ Therefore as right $A^e$-modules
$\frac{{}^1\!A^{\xi_f}\ot A^{op}}{D(I)({}^1\!A^{\xi_f}\ot
A^{op})}\cong {}^1\!A^{\xi_f}$. The proof is then complete. \qed

\begin{thm} \label{thm3} Let $\g$ be a finite dimensional Lie algebra. Then for any 2-cocycle $f\in Z^2(\g,\k)$, the following statements are equivalent.
\begin{itemize}
\item [(i)] The Sridharan
enveloping algebra $U_f(\g)$ is CY of dimension $d$.
\item [(ii)] The universal enveloping algebra $U(\g)$ is CY of dimension
  $d$.
\item [(iii)] $\dim\g=d$ and for any $x\in\g$, $\text{\rm tr}(\text{\rm ad}_\g(x))=0$.
\end{itemize}
\end{thm}
\proof Following Lemma \ref{lem5} it is sufficient to show that (i)
and (ii) are equivalent. We show first (ii) $\Rightarrow$ (i).
Assume that $U(\g)$ is CY of dimension $d$. Then $\dim(\g)=d$. Note
that $U(\g)$ is a cocommutative Hopf algebra. By Theorem
\ref{prop4}, $\RHom_{U(\g)}(\k,U(\g))\cong\k[d]$ as objects in the
derived category of complexes of right $U(\g)$-modules, where  $\k$
is the trivial right $U(\g)$-module. Once again we write $A$ for
$U_f(\g)$. Recall that $A^e$ is a free $U(\g)$-module. Now let
$P^\bullet$ be the Chevalley-Eilenburg resolution of the trivial
left $U(\g)$-module $\k$. Then $A^e\ot_{U(\g)} P^\bullet$ is the
standard resolution of $A$ as a left $A^e$-module (also see
\cite[Prop. 3]{Ka}). It follows that we have the following
isomorphisms in the derived category $D^\circ(A^e)$ of complexes of
right $A^e$-modules:
$$\begin{array}{ccl}
   \RHom_{A^e}(A,A^e)&\cong&\Hom_{A^e}(A^e\ot_{U(\g)}P^\bullet,A^e)\\
   &\cong&\Hom_{U(\g)}(P^\bullet,A^e)\\
&\cong&\Hom_{U(\g)}(P^\bullet,U(\g))\ot_{U(\g)} A^e\\
&\cong&\RHom_{U(\g)}(\k,U(\g))\ot_{U(\g)} A^e\\
&\cong&\k[-d]\ot_{U(\g)} A^e\\
&\overset{(a)}\cong&A[-d],
  \end{array}
$$
where the isomorphism $(a)$ follows from the right version of the
proof of \cite[Prop. 5.3]{Sr}. Therefore $U_f(\g)=A$ is a CY algebra
of dimension $d$.

(i) $\Rightarrow$ (ii). Assume that $A=U_f(\g)$ is CY of dimension
$d$. The first four isomorphisms above are still valid, and thus we
have
\begin{equation}\label{tag22}
\RHom_{A^e}(A,A^e)\cong\RHom_{U(\g)}(\k,U(\g))\ot_{U(\g)} A^e
\end{equation}
as right $A^e$-modules. Since $A$ is CY of dimension $d$ and $A^e$
is a free left $U(\g)$-module, $H^i\RHom_{U(\g)}(\k,U(\g))=0$ for
$i\neq d$ and $H^d\RHom_{U(\g)}(\k,U(\g))\neq0$. Hence
$\RHom_{U(\g)}(\k,U(\g))\cong \k^\xi[-d]$, where $\xi$ is, by
\cite[Prop. 6.3]{BZ}, the algebra automorphism $\xi:U(\g)\to U(\g)$
defined by $\xi(x)=x+\text{tr}(\text{ad}_\g(x))$ for all $x\in\g$.
Let $h=\text{tr}(\text{ad}_\g(-)):\g\to \k$. Then $h\in Z^1(\g,\k)$.
Now the 1-cocycle $h$ also defines an algebra automorphism
$\xi_f:A\to A$. Combining the isomorphisms in (\ref{tag22}) and the
isomorphism in Lemma \ref{lem22}, we  obtain the following
isomorphisms:
$$A[-d]\cong \k^\xi[-d]\ot_{U(\g)} A^e\cong {}^1\!A^{\xi_f}.$$
Thus we have an $A$-$A$-bimodule isomorphism: $A\cong {}^1\!A^{\xi_f}$. It follows that the isomorphism $\xi_f:A\to A$ must be
inner. That is, $\xi_f(a)=u^{-1}au$ for some unit $u\in A$. It is
easy to see that $u\in\k$. Therefore $\xi_f=id$ and
$h=\text{tr}(\text{ad}_\g(-))=0$. By Lemma \ref{lem5}, $U(\g)$
is CY. Moreover $U(\g)$ is of dimension $d$. \qed

The proof of Theorem \ref{thm3} yields a more general fact about
rigid dualizing complexes (for the definition, see \cite{VdB}) of
Sridharan enveloping algebras. Let $\g$ be a finite dimensional Lie
algebra, and $f\in Z^2(\g,\k)$ a 2-cocycle. Let $A=U_f(\g)$ as
before. By \cite[Cor. 8.7]{VdB} or \cite[Prop. 1.1]{Y}, the rigid
dualizing complex $R$ of $A$ exists. Moreover, $R$ is invertible and
$R^{-1}=\RHom_{A^e}(A,A^e)$. Notice that the linear map
$h=\text{tr}(\text{ad}_\g(-)):\g\to \k$ is a 1-cocycle of $\g$. As
early pointed out at the beginning of this section, $h$ defines both
an isomorphism $\xi$ on $U(\g)$ and an isomorphism $\xi_f$ on
$U_f(\g)$. Now we have the following corollary which generalizes
\cite[Theorem A]{Y} to Sridharan enveloping algebras.

\begin{cor} Let $\g$ be a Lie algebra of dimension $d$, and $f\in
Z^2(\g,\k)$ a 2-cocycle. Then the rigid dualizing complex of the
Sridharan enveloping algebra $U_f(\g)$ is ${}^1\!U_f(\g)^{\zeta_f}[d]$,
where $\zeta_f:\, U_f(\g)\to U_f(\g)$ is an algebra automorphism, and is defined by
$$\zeta_f(x)=x-\text{\rm tr}(\text{\rm ad}_\g(x)),\qquad \text{for all $x\in \g$}.$$
\end{cor}
\proof Let $A=U_f(\g)$. By \cite[Prop. 6.3]{BZ},
$\RHom_{U(\g)}(\k,U(\g))\cong\k^\xi[-d]$. Following the isomorphisms
in (\ref{tag22}), we have
$$\RHom_{A^e}(A,A^e)\cong \k^\xi[-d]\ot_{U(\g)}A^e.$$ By Lemma
\ref{lem22}, $\RHom_{A^e}(A,A^e)\cong {}^1\!A^{\xi_f}[-d]$.
Therefore the rigid dualizing complex of $A$ is
$R={}^1\!A^{\xi^{-1}_f}[d]$. Write $\zeta_f$ for $\xi^{-1}_f$, We
obtain the desired result. \qed

Now we focus on 3-dimensional CY Sridharan enveloping algebras.
By Theorem \ref{thm3}, such an algebra must be constructed from a
3-dimensional Lie algebra. Combining Proposition \ref{prop5},
Theorem \ref{thm3} and \cite[Theorem 1.3]{N0}, we may list all
possible nonisomorphic 3-dimensional CY Sridharan enveloping
algebras.

\begin{thm}\label{thm4} Let $U_f(\g)$ be a Sridharan enveloping algebra of a finite dimensional Lie algebra
$\g$. Then $U_f(\g)$ is CY of dimension 3 if and only if $U_f(\g)$
is isomorphic to $\k\langle x,y,z\rangle/(R)$ with the commuting
relations $R$ listed in the following table:\\
\center{
\begin{tabular}{c|ccc}\hline
{\rm Case}&$\quad \{x,y\}\quad$&$\quad\{x,z\}\quad$&$\quad\{y,z\}\quad$ \\
\hline $1$&$z$&$-2x$&$2y$\\
$2$&$y$&$-z$&$0$\\
$3$&$z$&$0$&$0$ \\
$4$&$0$&$0$&$0$ \\
$5$&$y$&$-z$&$1$\\
$6$&$z$&$1$&$0$\\
$7$&$1$&$0$&$0$\\
\hline
\end{tabular}}

where $\{x,y\}=x y-y x$.
\end{thm}

Note that in the above table the cases 1--4 give the CY
universal enveloping algebras listed in Proposition \ref{prop5}.

Let $\mathfrak{g}$ be a finite dimensional Lie algebra, and $f\in
Z^2(\g,\k)$ a 2-cocycle. The Sridharan enveloping algebra
$U_f(\mathfrak{g})$ is a PBW-deformation of the polynomial algebra
$\k[x_1,\dots,x_n]$ where $n=\dim\g$ (cf. \cite{Sr,BT,PP}).
Conversely, a PBW-deformation of a polynomial algebra is exactly a
Sridharan enveloping algebra (cf. \cite{Sr,N}). It is shown in
\cite[Theorem 3.6]{BT} that if a PBW-deformation of a 3-dimensional
graded CY algebra is defined by a potential,  then the deformed
algebra is also CY of dimension 3. Whether the converse is true or
not is not shown in \cite{BT}. However, a CY Sridharan enveloping
algebra $U_f(\mathfrak{g})$ of a 3-dimensional Lie algebra is always
defined by a potential. In fact, all 3-dimensional Lie algebras
whose Sridharan enveloping algebras are CY are listed in Theorem
\ref{thm4}. It is easy to check that the defining relations of these
algebras satisfy the condition in \cite[Theorem 3.2]{BT}. Hence any
3-dimensional CY Sridharan enveloping algebra $A$ is defined by a
potential. That is; $A\cong\k \langle
x,y,z\rangle/(\frac{\partial\Phi}{\partial_x},\frac{\partial\Phi}{\partial_y},\frac{\partial\Phi}{\partial_z})$,
where $\Phi\in \k \langle x,y,z\rangle/[\k \langle x,y,z\rangle,\k
\langle x,y,z\rangle]$ is a potential.

In fact, we can write down the potentials corresponding to the
Sridharan enveloping algebras of the Lie algebras in Theorem
\ref{thm4} respectively:

(1) $\Phi=xyz-yxz-\frac{1}{2}z^2-2xy$;

(2) $\Phi=xyz-yxz-yz$;

(3) $\Phi=xyz-yxz-\frac{1}{2}z^2$;

(4) $\Phi=xyz-yxz$;

(5) $\Phi=xyz-yxz-yz-x$;

(6) $\Phi=xyz-yxz-\frac{1}{2}z^2-y$;

(7) $\Phi=xyz-yxz-z$.

Hence a PBW-deformation $A$ of the polynomial algebra $\k[x,y,z]$ is
CY if and only if $A$ is defined by a potential. This phenomenon
does not occur accidentally. Travis Schedler shows in
\cite{Sc} that any 3-dimensional CY PBW-deformation of a
3-dimensional graded CY algebra (associated to a finite quiver) must
be defined by a potential. Combining with the results of \cite{BT},
we then obtain that a PBW-deformation $A$ of a 3-dimensional graded
CY algebra (associated to a finite quiver) is CY of dimension 3 if
and only if $A$ is defined by a potential.

\vspace{5mm}

\subsection*{Acknowledgement} The authors would like to thank Raf
Bocklandt and Travis Schedler for useful conversations. In
particular, the first author thanks Travis Schedler for sharing his
ideas with him and showing him the manuscript \cite{Sc}. The work is
supported by an FWO-grant and NSFC (No. 10801099).

\vspace{5mm}

\bibliography{}

\end{document}